\documentclass[12pt]{article}
\usepackage{setspace}
\usepackage[margin=1in]{geometry}
\usepackage{amsmath}
\usepackage{amssymb}
\usepackage{graphicx,multirow}
\usepackage{algorithm}
\usepackage{algorithmic}
\usepackage{subfigure}
\usepackage{multicol,float}
\usepackage{bbm}
\usepackage{url}
\usepackage{tikz,pgfplots}
\usepackage{amsthm}
\usepackage{authblk}
\usepackage[round,sort]{natbib}
\usepackage[pdfborder={0 0 0}]{hyperref}
\newcommand{\ind}{\mathbbm{1}}

\def\RR{\mathbb{R}}

\def\EE{\mathbb{E}}
\def\PP{\mathbb{P}}

\DeclareMathOperator*{\median}{median}

\begin{document}

\title{Mini-batch stochastic gradient descent\\with dynamic sample sizes}

\author{Michael R. Metel}

\affil{Laboratoire de Recherche en Informatique, Universit\'e Paris-Sud, Orsay, France\\
\url{metel@lri.fr}}

\maketitle

\begin{abstract}
We focus on solving constrained convex optimization problems using mini-batch stochastic gradient descent. Dynamic sample size rules are presented which ensure a descent direction with high probability. Empirical results from two applications show superior convergence compared to fixed sample implementations.
\end{abstract}

\section{Introduction}

\doublespace

We are interested in the following optimization problem,

$$\min\limits_{x\in X} f(x):=\EE_y\left[F(x,y)\right]$$

where $X\subset \RR^m$ is a convex feasible region, $y\in \RR^n$ is a random vector following a distribution from which we can generate i.i.d. samples, and $F(\cdot,y)$ is convex. If $f(\cdot)$ is finite valued in a neighbourhood of $x$, then
$$g(x):=\nabla f(x)=\EE_y\left[\nabla F(x,y)\right],$$

which we will assume throughout, see \cite{shap2009}.\\

We can solve this problem using a gradient descent algorithm, with updates approximating

$$x^{i+1}=\pi(x^i-\eta^i g(x^i)),$$

where $\pi(\cdot)$ is the Euclidian projection onto $X$, and $\eta^i$ is a chosen step size or learning rate. This formulation is not practical for large scale problems, requiring multidimensional integration each iteration. Stochastic gradient descent (SGD) algorithms take a sample, $y^i$, and use $\nabla F(x^i,y^i)$ in the iterative scheme

$$x^{i+1}=\pi(x^i-\eta^i\nabla F(x^i,y^i)).$$

This type of algorithm dates back to \cite{robbins1951}, where convergence was established for learning rates of the form $\eta^i=\Theta(\frac{1}{i})$. Since then, there have been many improvements and new techniques employed to improve performance of this iterative algorithm, such as momentum \citep{qian1999}, averaging \citep{nem2009}, and more recently a number of adaptive learning rate methods, such as Adam \citep{kingma2014}.\\

In expectation, $\EE\left[\nabla F(x^i,y^i)\right]=g(x^i)$, but we can imagine using a single sample will be noisy, with better estimates following from using a mini-batch of $N^i$ instances of $\nabla F(x^i,y^i)$,

\begin{equation} \label{eq:MB}
x^{i+1}=\pi(x^i-\frac{\eta^i}{N^i}\sum_{j=1}^{N^i} \nabla F(x^i,y_j^i)).
\end{equation}

We can also expect there to be a diminishing return on the sample size, and that at a certain level the computational cost of using more samples will be greater than the improved performance from using a more accurate gradient estimate. There has been some research examining what sample size to use, and in particular variable sample sizes which prove convergence using a geometrically increasing sample size, see \cite{hash2014} and \cite{byrd2012}. In addition, \cite{byrd2012} presented a condition which ensures that the estimated gradient is a descent direction, with a heuristic to approximate the appropriate sample size to satisfy it. In this work we are interested in developing a sample size rule which uses as little samples as possible while ensuring a descent direction with high probability, to achieve improved convergence in practice.

\section{Dynamic sample sizes}

Let $$\hat{g}(x^i)=\frac{1}{N^i}\sum_{j=1}^{N^i}\nabla F(x^i,y_j^i)$$

be our estimate of the gradient, from which we can calculate its sample covariance matrix,

$$\hat{\Sigma}(x^i)=\frac{1}{N^i(N^i-1)}\sum_{j=1}^{N^i}(\nabla F(x^i,y_j^i)-\hat{g}(x^i))(\nabla F(x^i,y_j^i)-\hat{g}(x^i))^T\nonumber$$

When weighing the trade-off between gradient estimate accuracy and computational cost, a base level of accuracy is moving in a descent direction. By the central limit theorem, we know that $\hat{g}(x^i)-g(x^i)\sim N(0,\Sigma(x^i))$ approximately, where $\Sigma(x^i)=\frac{1}{N^i}\text{Cov}(\nabla F(x^i,y))$. Replacing the actual covariance of $\hat{g}(x^i)$ with its sample estimate, we can estimate confidence intervals of $g(x^i)$ with distribution

$$g(x^i)\sim t_{N^i-1}(\hat{g}(x^i),\hat{\Sigma}(x^i)).$$

Further, we assume $N^i$ is large enough so as to make the change to a normal approximation insignificant for our purposes, and model

$$g(x^i)\sim N(\hat{g}(x^i),\hat{\Sigma}(x^i))$$

\subsection{Per dimension update}

Let us consider the $k^{th}$ entry of $\hat{g}(x^i)$, and assume that it is greater than zero, $\hat{g}(x^i)_k>0$. We can then estimate the probability that $g(x^i)_k$ is as well,

$$\PP\left(g(x^i)_k>0\right)=\PP\left(\frac{g(x^i)_k-\hat{g}(x^i)_k}{\sqrt{\hat{\Sigma}(x^i)_{kk}}}>\frac{-\hat{g}(x^i)_k}{\sqrt{\hat{\Sigma}(x^i)_{kk}}}\right)
=1-\Phi\left(\frac{-\hat{g}(x^i)_k}{\sqrt{\hat{\Sigma}(x^i)_{kk}}}\right)=\Phi\left(\frac{\hat{g}(x^i)_k}{\sqrt{\hat{\Sigma}(x^i)_{kk}}}\right).$$

Likewise, for $\hat{g}(x^i)_k\leq 0$,

$$\PP\left(g(x^i)_k\leq 0\right)=\PP\left(\frac{g(x^i)_k-\hat{g}(x^i)_k}{\sqrt{\hat{\Sigma}(x^i)_{kk}}}\leq\frac{-\hat{g}(x^i)_k}{\sqrt{\hat{\Sigma}(x^i)_{kk}}}\right)
=\Phi\left(\frac{-\hat{g}(x^i)_k}{\sqrt{\hat{\Sigma}(x^i)_{kk}}}\right).$$

So in general the probability of moving in the direction of descent by moving in the direction $-\hat{g}(x^i)_k$ equals $\Phi\left(\frac{|\hat{g}(x^i)_k|}{\sqrt{\hat{\Sigma}(x^i)_{kk}}}\right)$.\\

In this subsection we consider using a different sample size, $N^i_k$, for each partial derivative, and desire a movement in a descent direction with probability $1-\alpha$ for $\alpha\in (0,0.5)$ in each dimension. Assume using a current number $N_k^i$ of samples, $\Phi\left(\frac{|\hat{g}(x^i)_k|}{\sqrt{\hat{\Sigma}(x^i)_{kk}}}\right)<1-\alpha$.
We can achieve a higher probability direction by decreasing $\sqrt{\hat{\Sigma}(x^i)_{kk}}\approx \sqrt{\frac{1}{N_k^i}\text{Cov}(\nabla F(x^i,y))}$ by increasing $N^i_k$. For a decreased standard deviation, $\sqrt{\hat{\Sigma}(x^{i'})_{kk}}$, such that
$\frac{\sqrt{\hat{\Sigma}(x^{i'})_{kk}}}{\sqrt{\hat{\Sigma}(x^i)_{kk}}}\leq \theta$ for $0<\theta<1$, we must choose an increased sample size $N_k^{i'}$ such that $\frac{N_k^i}{N_k^{i'}}\leq \theta^2$. In particular, if we want $\Phi\left(\frac{|\hat{g}(x^i)_k|}{\sqrt{\hat{\Sigma}(x^{i'})_{kk}}}\right)\geq 1-\alpha$, or
$\sqrt{\hat{\Sigma}(x^{i'})_{kk}}\leq \frac{|\hat{g}(x^i)_k|}{\Phi^{-1}(1-\alpha)}$, then
$\frac{\sqrt{\hat{\Sigma}(x^{i'})_{kk}}}{\sqrt{\hat{\Sigma}(x^i)_{kk}}}\leq\frac{|\hat{g}(x^i)_k|}{\sqrt{\hat{\Sigma}(x^i)_{kk}}\Phi^{-1}(1-\alpha)}$, implying
$\frac{N^i_k}{N^{i'}_k}\leq\left(\frac{|\hat{g}(x^i)_k|}{\sqrt{\hat{\Sigma}(x^i)_{kk}}\Phi^{-1}(1-\alpha)}\right)^2$,
and so we choose
\begin{equation} \label{eq:PDUR}
N_k^{i+1}=\left\lceil N_k^{i}\frac{\hat{\Sigma}(x^i)_{kk}\left(\Phi^{-1}(1-\alpha)\right)^2}{\hat{g}(x^i)^2_k}\right\rceil
\end{equation}

as the sample size in the next iteration. As we want to use samples sparingly, in the case where $\Phi\left(\frac{|\hat{g}(x^i)_k|}{\sqrt{\hat{\Sigma}(x^i)_{kk}}}\right)>1-\alpha$, we also use (\ref{eq:PDUR}) to decrease our sample size for the next iteration. Implicit in this scheme is the assumption that the variance in gradient samples do not vary significantly from iteration to iteration, making our estimate of $N^{i+1}$ using information from iteration $i$ accurate.\\

In empirical testing we generated estimates $\hat{g}(x^{i+1})_k$ using separate $N^{i+1}_k$ in each dimension, as well as using $\max\limits_j N^{i+1}_j$ samples for all dimensions, but found in practice using $\median\limits_j N^{i+1}_j$ to be most effective, given its simplicity and ability to avoid at times large sample size outliers.

\subsection{Single update}

A less stringent approach to sample size selection is to require that $-\hat{g}(x^i)$ is a direction of descent with high probability in aggregate, which we know holds if $\hat{g}(x^i)^Tg(x^i)>0$. From our assumptions,

$$\hat{g}(x^i)^Tg(x^i)\sim N(\hat{g}(x^i)^T\hat{g}(x^i),\hat{g}(x^i)^T\hat{\Sigma}(x^i)\hat{g}(x^i))$$

Following the same steps as in the previous subsection, we get the update rule

$$N^{i+1}=\left\lceil N^{i}\frac{\hat{g}(x^i)^T\hat{\Sigma}(x^i)\hat{g}(x^i)(\Phi^{-1}(1-\alpha))^2}{(\hat{g}(x^i)^T\hat{g}(x^i))^2}\label{eq:SUR}\right\rceil.$$

In empirical testing we found improved performance by assuming estimated gradient terms are independent, leaving us with only a diagonal covariance matrix to estimate.

\section{Numerical experiments}
\label{sec:ES}

We compare the performance of the dynamic sampling approaches to fixed sample implementations of the basic mini-batch SGD algorithm (\ref{eq:MB}) with a decreasing learning rate of $\eta^i=\frac{1}{i}$ and Adam using the default parameter values presented in \citep{kingma2014}, namely $\eta^i = 0.001, \beta_1 = 0.9, \beta_2 = 0.999, \text{ and } \epsilon = 10^{-8}$, in the following two applications.

\subsection{Newsvendor problem}

We first consider a single period multi-product newsvendor problem with an exponential utility function and correlated demand \citep{choi2011}. The decision variables are $x_j$, the amount of product we order, at a cost per unit $c_j$, with selling price $p_j$, and uncertain demand $D_j$. Given $n$ products our random profit is

$$\sum_{j=1}^np_j\min\{x_j,D_j\}-c_jx_j.$$

The exponential utility function is of the form $u(z)=-e^{-\lambda z}$ where $\lambda$ is a risk aversion parameter. The optimization problem is as follows.

\begin{alignat}{6}
\label{eq:NV}
&\max&&\text{ }\EE[-e^{-\lambda\left(\sum_{j=1}^np_j\min\{x_j,D_j\}-c_jx_j\right)}]\nonumber\\
&\mbox{s.t. }&&x_i\geq 0\nonumber
\end{alignat}

and

$$\nabla F(x,D)_i=\lambda(p_i\ind_{\{x_i<D_i\}}-c_i)e^{-\lambda\left(\sum_{j=1}^np_j\min\{x_j,D_j\}-c_jx_j\right)}.$$

We generated random data for 50 products, using values similar to \cite{choi2011}. The prices $p_i$ were uniformly sampled from $[15,30]$, $c_i=10$, and $\lambda=0.02$. The random demands follow a log-normal distribution generated from a normal distribution with $\mu_i=3$, $\sigma_i$ uniformly sampled from $[0.4724,1.2684]$ to achieve coefficients of variance between $[0.5,2]$, and constant correlations of $\rho_{i,j}=0.25$ between variables.

\subsection{Call and put options portfolio problem}

The second application is finding the optimal growth portfolio, see \cite{estrada2010}, of European call and put options with stock returns following geometric Brownian motions,

$$dS^j_t=\mu_j S^j_t dt +\sigma_j S^j_t dW^j_t$$

where $\mu_j$ is the expected stock return, $\sigma_j$ is the standard deviation of the stock return, and the Brownian motions have correlation

$$\EE(dW^j_tdW^k_t)=\rho_{j,k}dt.$$

We invest in at the money call and put options on each stock at time $t=0$ to maximize our return at time $t=1$, with random stock prices equal to $S^j_1=S^j_0e^{(\mu_j-\frac{1}{2}\sigma_j^2)+\sigma_jW^j_1}$, call option payoffs $C_1^j=\max(S_1^j-S_0^j,0)$, and put option payoffs
$P_1^j=\max(S_0^j-S_1^j,0)$.\\

The optimization problem is to maximize the expected logarithm of wealth. The decision variable $x^C_j$ is the fraction of wealth invested in $C^j$, $x^P_j$ is the fraction of wealth invested in $P^j$, and $r$ is the interest rate received from lending money.

\begin{alignat}{6}
&\max&&\text{ }\EE\log\left(1+r+\sum_{i=1}^m x_i^C\left(\frac{C_1^i}{C_0^i}-(1+r)\right)+
x_i^P\left(\frac{P_1^i}{P_0^i}-(1+r)\right)\right)\nonumber\\
&\mbox{s.t. }&&\sum_i^m x_i^C+x_i^P\leq 1\nonumber\\
&&&x_i^C,x_i^P\geq 0\nonumber
\end{alignat}

In implementing the SGD algorithm,

$$\nabla F(x,C_1,P_1)_i^C=\frac{\frac{C_1^i}{C_0^i}-(1+r)}{1+r+\sum_{i=1}^m x_i^C\left(\frac{C_1^i}{C_0^i}-(1+r)\right)+
x_i^P\left(\frac{P_1^i}{P_0^i}-(1+r)\right)}$$

$$\nabla F(x,C_1,P_1)_i^P=\frac{\frac{P_1^i}{P_0^i}-(1+r)}{1+r+\sum_{i=1}^m x_i^C\left(\frac{C_1^i}{C_0^i}-(1+r)\right)+
x_i^P\left(\frac{P_1^i}{P_0^i}-(1+r)\right)}$$

We simulated a universe of 50 stocks, using the $ev-\hat{e}\hat{v}$ methodology of \cite{hirschberger2007} to generate two random covariance matrices, with the parameterization provided from randomly selected stocks from the S\&P SuperComposite 1500. We simulated our estimate of the covariance matrix of stock returns, $\Sigma$, and the market's, $\Sigma^m$. Given $\Sigma^m$ and $r$, we calculated $C_0$ and $P_0$ assuming the market follows the Black-Scholes model of option pricing. With $\Sigma$, we simulated our estimate of expected stock returns $\mu_{i}$ with magnitude uniformly distributed between $[0,2]*\sigma_{ii}$, and positive with probability $0.75$.

\section{Results}
\label{sec:R}

All experiments were done on a Windows 10 Home 64-bit, Intel Core i5-7200U 2.5GHz processor with 8 GB of RAM, in Matlab R2017a. Mini-batch sample sizes used generally vary between 32 and 512 samples \citep{keskar2016}, so both experiments consisted of testing the per dimension and single update rules against fixed sample sizes of 32, 256, and 512. We found superior convergence using the dynamic sampling approaches in both applications. Below are plots of the objective value through time.

 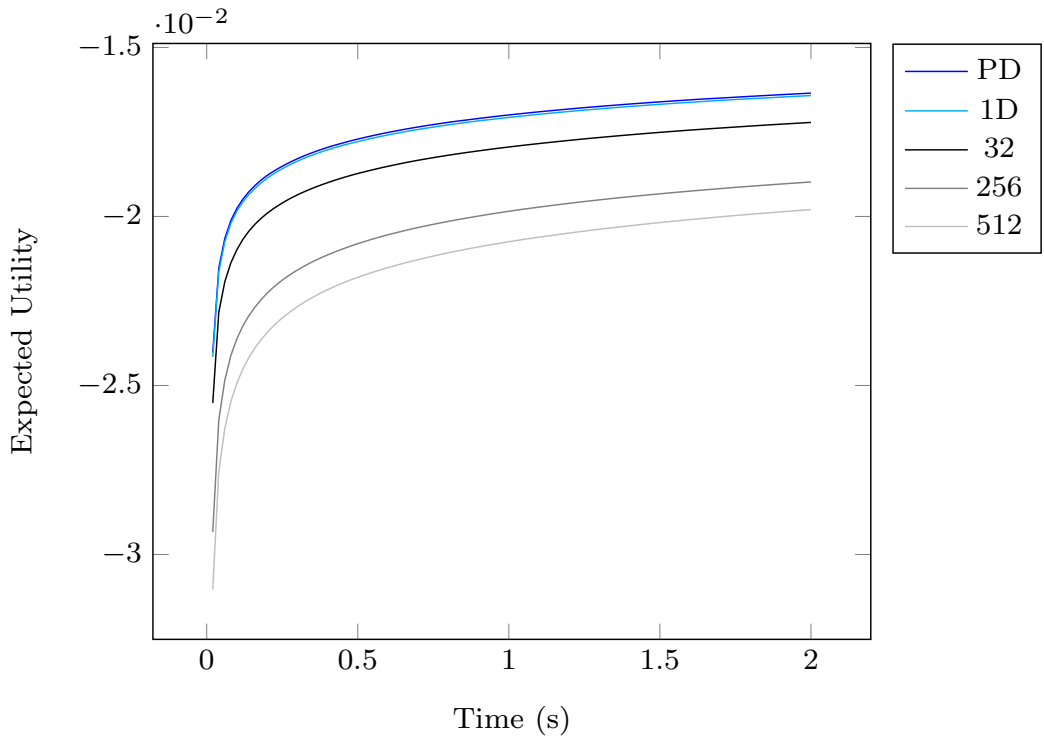
\begin{figure}[H]
 \centerline{
\resizebox{0.85\textwidth}{!}{
\begin{tikzpicture}
\scriptsize
    \begin{axis}[xlabel=Time (s),ylabel=Expected Utility,
                legend style={legend pos=outer north east}]
                \addplot[mark=none,mark size=0.75,draw=blue]
                table[x=T,y=BPD]
            {resultsNB4.dat};
            \addplot[mark=none,mark size=0.75,draw=cyan]
                table[x=T,y=B1D]
            {resultsNB4.dat};
            \addplot[mark=none,mark size=0.75,draw=black]
                table[x=T,y=B32]
            {resultsNB4.dat};
            \addplot[mark=none,mark size=0.75,draw=gray]
                table[x=T,y=B256]
            {resultsNB4.dat};
            \addplot[mark=none,mark size=0.75,draw=lightgray]
                table[x=T,y=B512]
            {resultsNB4.dat};
    \legend{PD,1D,32,256,512}
    \end{axis}
\end{tikzpicture}}}
\caption{Newsvendor problem using basic SGD} \label{T3}
\end{figure}
 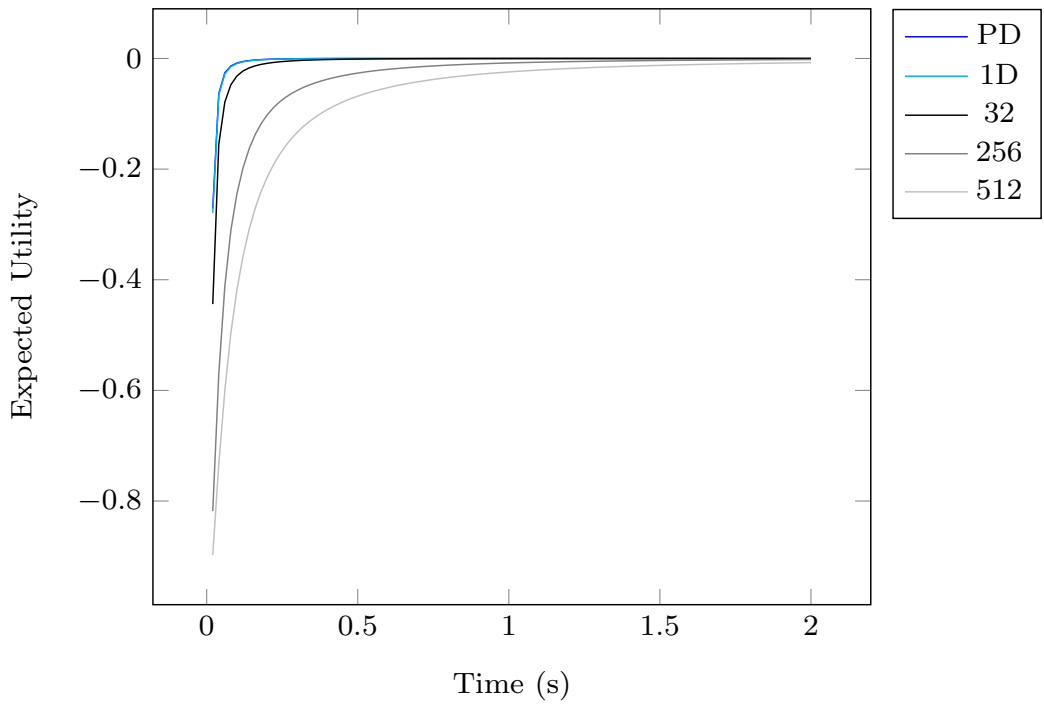
\begin{figure}[H]
 \centerline{
\resizebox{0.85\textwidth}{!}{
\begin{tikzpicture}
\scriptsize
    \begin{axis}[xlabel=Time (s),ylabel=Expected Utility,
                legend style={legend pos=outer north east}]
                \addplot[mark=none,mark size=0.75,draw=blue]
                table[x=T,y=APD]
            {resultsNB4.dat};
            \addplot[mark=none,mark size=0.75,draw=cyan]
                table[x=T,y=A1D]
            {resultsNB4.dat};
            \addplot[mark=none,mark size=0.75,draw=black]
                table[x=T,y=A32]
            {resultsNB4.dat};
            \addplot[mark=none,mark size=0.75,draw=gray]
                table[x=T,y=A256]
            {resultsNB4.dat};
            \addplot[mark=none,mark size=0.75,draw=lightgray]
                table[x=T,y=A512]
            {resultsNB4.dat};
    \legend{PD,1D,32,256,512}
    \end{axis}
\end{tikzpicture}}}
\caption{Newsvendor problem using Adam} \label{T4}
\end{figure}
 \begin{figure}[H]
 \centerline{
\resizebox{0.85\textwidth}{!}{
\begin{tikzpicture}
\scriptsize
    \begin{axis}[xlabel=Time (s),ylabel=Expected Utility,
                legend style={legend pos=outer north east}]
                \addplot[mark=none,mark size=0.75,draw=blue]
                table[x=T,y=BPD]
            {resultsO4.dat};
            \addplot[mark=none,mark size=0.75,draw=cyan]
                table[x=T,y=B1D]
            {resultsO4.dat};
            \addplot[mark=none,mark size=0.75,draw=black]
                table[x=T,y=B32]
            {resultsO4.dat};
            \addplot[mark=none,mark size=0.75,draw=gray]
                table[x=T,y=B256]
            {resultsO4.dat};
            \addplot[mark=none,mark size=0.75,draw=lightgray]
                table[x=T,y=B512]
            {resultsO4.dat};
    \legend{PD,1D,32,256,512}
    \end{axis}
\end{tikzpicture}}}
\caption{Options portfolio using basic SGD} \label{T1}
\end{figure}
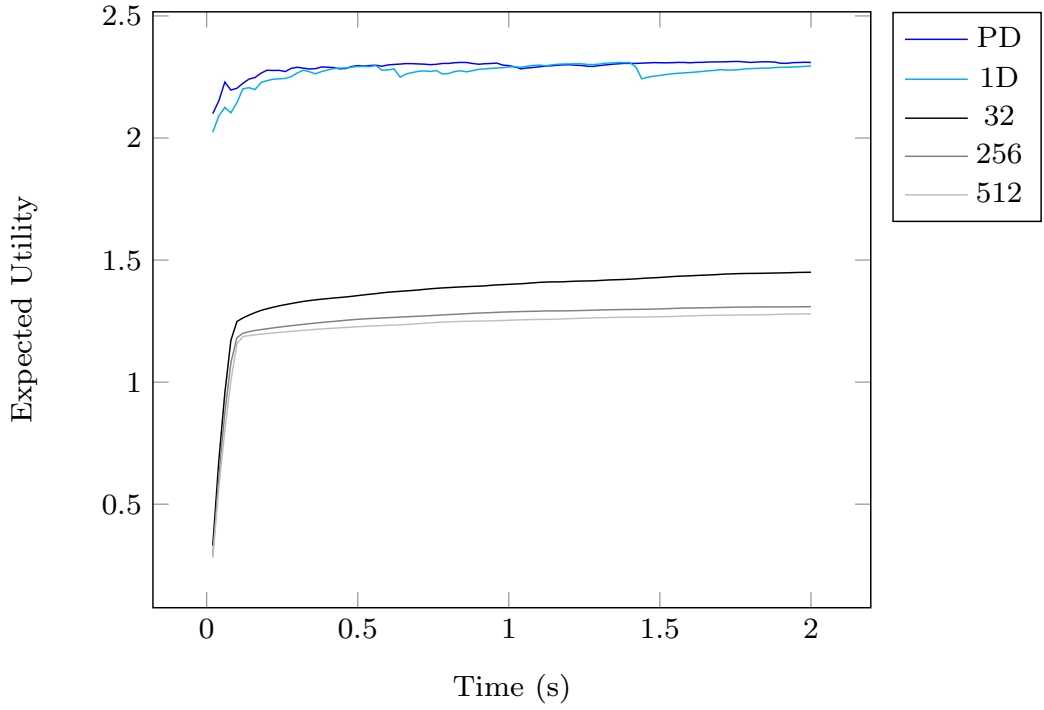
 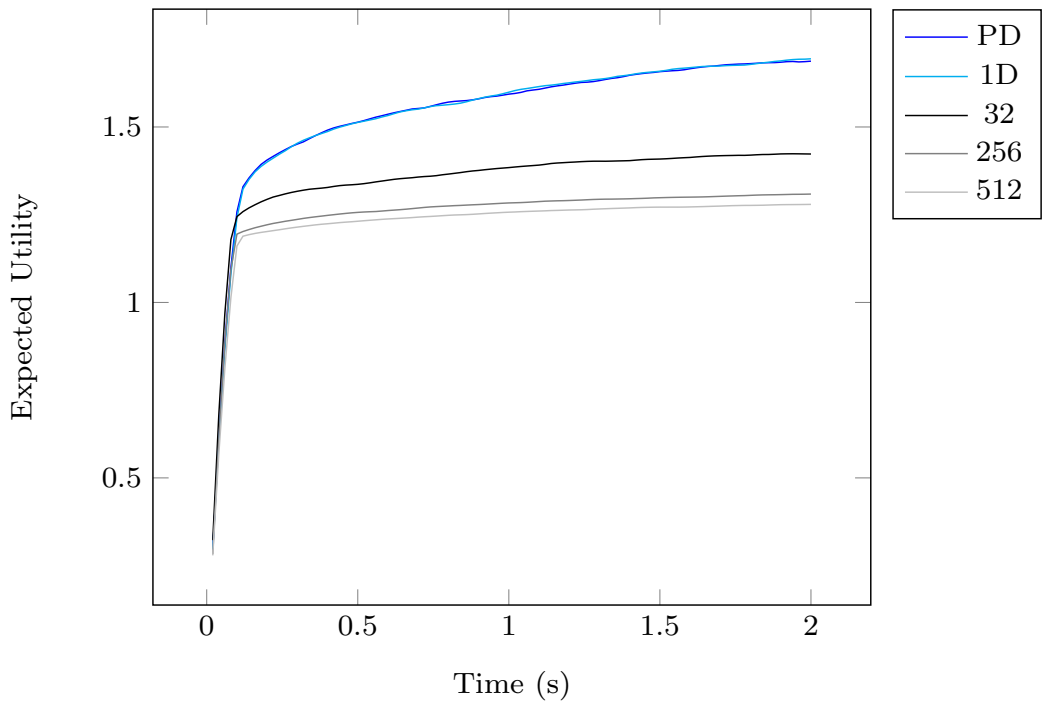
\begin{figure}[H]
 \centerline{
\resizebox{0.85\textwidth}{!}{
\begin{tikzpicture}
\scriptsize
    \begin{axis}[xlabel=Time (s),ylabel=Expected Utility,
                legend style={legend pos=outer north east}]
                \addplot[mark=none,mark size=0.75,draw=blue]
                table[x=T,y=APD]
            {resultsO4.dat};
            \addplot[mark=none,mark size=0.75,draw=cyan]
                table[x=T,y=A1D]
            {resultsO4.dat};
            \addplot[mark=none,mark size=0.75,draw=black]
                table[x=T,y=A32]
            {resultsO4.dat};
            \addplot[mark=none,mark size=0.75,draw=gray]
                table[x=T,y=A256]
            {resultsO4.dat};
            \addplot[mark=none,mark size=0.75,draw=lightgray]
                table[x=T,y=A512]
            {resultsO4.dat};
    \legend{PD,1D,32,256,512}
    \end{axis}
\end{tikzpicture}}}
\caption{Options portfolio using Adam} \label{T2}
\end{figure}

\section{Conclusion}
\label{s:C}
We have presented two rules to dynamically select sample sizes in SGD algorithms to ensure a direction of descent with high probability using as little samples as possible. Superior convergence was found compared to fixed sample approaches with the per dimension update rule having superior performance overall in two test applications.

\bibliographystyle{plainnat}
\bibliography{StoGrad}

\end{document}